\documentclass[11pt,final,a4paper]{amsart}
\usepackage{amssymb}
\usepackage{amsmath}
\usepackage{amsfonts}

\newcommand{\F}{\mathbb{F}}

\newtheorem{main-dummy}{Main-Dummy}
\newtheorem{dummy}{Dummy}

\newtheorem{main-theorem}[main-dummy]{Theorem}

\newtheorem{lemma}[dummy]{Lemma}
\newtheorem*{lemma*}{Lemma}
\newtheorem{theorem}[dummy]{Theorem}

\newtheorem{cor}[dummy]{Corollary}

\theoremstyle{definition}

\newtheorem{example}[dummy]{Example}

\theoremstyle{remark}

\begin{document}
\bibliographystyle{amsalpha}
\author{Sandro Mattarei}

\email{mattarei@science.unitn.it}

\urladdr{http://www-math.science.unitn.it/\~{ }mattarei/}

\address{Dipartimento di Matematica\\
  Universit\`a degli Studi di Trento\\
  via Sommarive 14\\
  I-38050 Povo (Trento)\\
  Italy}

\title[Root multiplicities and  number of nonzero coefficients]
{Root multiplicities and  number of nonzero coefficients of a polynomial}

\begin{abstract}
It is known that the weight (that is, the number of nonzero coefficients) of a univariate polynomial
over a field of characteristic zero
is larger than the multiplicity of any of its nonzero roots.
We extend this result to an appropriate statement in positive characteristic.
Furthermore, we present a new proof of the original result, which produces also
the exact number of monic polynomials of a given degree for which the bound is attained.
A similar argument allows us to determine the number of monic polynomials
of a given degree, multiplicity of a given nonzero root, and number of nonzero coefficients,
over a finite field of characteristic larger than the degree.
\end{abstract}

\date{January 20, 2006}



\thanks{Partially supported by MIUR-Italy via PRIN 2003018059
 ``Graded Lie algebras and pro-$p$-groups: representations,
 periodicity and derivations''.}

\maketitle

\thispagestyle{empty}

\section{Introduction}

Call {\em weight} of a univariate polynomial over a field the number of its nonzero coefficients.
This concept is relevant to various branches of mathematics.
To mention just one example, polynomials over finite fields are widely used in coding theory and cryptography,
and in certain situations those with low weight are preferable, because they allow faster computations.

There are various relationships between the coefficients of a polynomial and the distribution of its roots.
A very simple one involving the weight is the following, which is Lemma~1 of~\cite{Bri}.

\begin{lemma}\label{lemma:weight}
The weight of a polynomial $f(x)$ over a field of characteristic zero exceeds at least by one
the multiplicity of any of its nonzero roots.
\end{lemma}

The simplicity of this result suggests that it may be folklore.
For example, an analogue in positive characteristic appears
in the paper~\cite{H-BK} cited below.
Brindza's proof is by induction on the multiplicity, and coincides in essence with the first two paragraphs of the proof
of our Theorem~\ref{thm:HBK_extended} given in the next section.
The polynomial $(x-1)^n$ shows that the lower bound for the weight of $f$ given in Lemma~\ref{lemma:weight} is best possible.

An immediate consequence of Lemma~\ref{lemma:weight}, as drawn in Theorem~1 of~\cite{Bri},
is that the weight of a polynomial $f$ with $f(0)\neq 0$ exceeds at least by one the ratio of the degree of $f$
over the number of its distinct roots in the field of definition.
Note that the assumption that $0$ is not a root of $f$ is harmless and natural, because the weight of $f$
remains unaffected after dividing $f$ by a power of $x$.
Again, the polynomial $(x-1)^n$ shows that this result is sharp, as noted in~\cite{Bri}.
We mention here the more general example $(x^d-1)^n$, which has each of its $d$ distinct roots
(in a splitting field) with multiplicity $n$,
and has weight $n+1=nd/d+1$.

Lemma~\ref{lemma:weight} remains true in prime characteristic $p$
provided the degree of $f$ is less than $p$,
and in this form it is proved and used by Heath-Brown and Konyagin in~\cite{H-BK}.
In fact, it is enough to assume that the root multiplicity involved is less than $p$,
as stated and applied in~\cite[Lemma~9]{Mat:binomial-recurrent}.
The polynomial $x^p-1$ shows that this assumption cannot be further weakened.
The first goal of this note is to prove the best possible conclusion when this assumption is removed.

\begin{theorem}\label{thm:HBK_extended}
Let $f(x)$ be a polynomial over a field of characteristic $p$,
having a nonzero root $\xi$ with (positive) multiplicity exactly $k$.
Then $f(x)$ has weight at least $\prod_t(k_t+1)$, where
$k=\sum_tk_tp^t$ is the $p$-adic expansion of $k$.
\end{theorem}

The conclusion of Theorem~\ref{thm:HBK_extended} is best possible, because $\prod_t(k_t+1)$ is the weight of the polynomial
$(x-1)^k=\prod_t(x^{p^t}-1)^{k_t}$.

Our second goal is to produce another proof of Lemma~\ref{lemma:weight}
based on a different idea from that of~\cite{Bri}.
Instead of working by induction on $k$, we deduce the result by linear algebra directly from the fact that $(x-1)^k$ has weight $k+1$.
Although this approach does not seem to easily extend to a proof of Theorem~\ref{thm:HBK_extended},
it has other advantages, because it can be used to produce additional information.
To illustrate this claim, we extend the argument
to compute the number of polynomials of degree $n$ (in characteristic zero), having a nonzero root
(which may conveniently be taken to be $1$ or $-1$, see the first paragraph in Section~\ref{sec:proofs})
of multiplicity $k$ and minimal weight $k+1$.

\begin{theorem}\label{thm:count_extremal}
Over a field of characteristic zero there are exactly $\binom{n}{k}$ monic polynomials $f(x)$ of degree $n$
which are multiples of $(x+1)^k$ and have weight $k+1$.
Each of them is determined uniquely by the set of degrees of its nonzero monomials.
\end{theorem}

Theorem~\ref{thm:count_extremal} remains true over a field of characteristic $p>n$, and becomes then a special case
of Theorem~\ref{thm:count} below.
We illustrate some of the difficulties which one runs into without this assumption
by working out the case $n=p$ in Example~\ref{ex:n=p}.

The argument used to prove Theorem~\ref{thm:count_extremal}
can be extended to give a description of the set of monic polynomials $f(x)$ of degree $n$
which are multiples of $(x+1)^k$ and have given weight $w$, for any $w>k+1$.
While such a set is infinite in characteristic zero,
over a finite field we can compute its finite cardinality explicitly,
again under the additional assumption that $n<p$.
This is an immediate consequence of the following more precise result.

\begin{theorem}\label{thm:count}
Let $\F_q$ be the field of $q$ elements, let $p$ be its characteristic,
and let $1\le k\le n<p$.
For all positive integers $w$, set
\[
M_{w}=\sum_{v\ge k}(-1)^{w-v-1}\binom{w-1}{v}q^{v-k}.
\]
In particular, $M_w$ vanishes for $w\le k$.

Let $J$ be a subset of $\{0,1,\ldots,n-1\}$ with $|J|=w-1$.
Then $M_w$ equals
the number of polynomials in $\F_q[x]$
of the form
\[
x^n+\sum_{j\in J}f_jx^j=f(x)=(x+1)^kg(x)
\]
for some $g(x)\in\F_q[x]$, with $f_j\neq 0$ for all $j\in J$.
\end{theorem}

\begin{cor}\label{cor:count}
Let $\F_q$ be the field of $q$ elements, let $p$ be its characteristic,
and let $1\le k\le n<p$.
Define $M_w$ as in Theorem~\ref{thm:count}.
Then the number of monic polynomials $f(x)\in\F_q[x]$ of degree $n$
which are multiples of $(x+1)^k$ and have weight exactly $w$, equals
$\binom{n}{w-1}M_w$.
\end{cor}

\section{Proofs}\label{sec:proofs}

If $\xi$ is a nonzero root of $f(x)$ with multiplicity exactly $k$, then $f(\xi^{-1}x)$
has $1$ as a root with the same multiplicity, and has the same weight as $f(x)$.
Hence we assume that $\xi=1$ in the following proof.

\begin{proof}[Proof of Theorem~\ref{thm:HBK_extended}]
We proceed by induction on $k$.
The case $k=1$ being obvious, assume that $k>1$.
By dividing $f(x)=\sum_i f_ix^{i}$ by a suitable power of $x$, which leaves its weight unchanged,
we may assume that $f_0\neq 0$.

If $k<p$
the derivative
\[
f'(x)=\sum_i if_ix^{i-1}
\]
has $1$ as a root with multiplicity exactly $k-1$,
and has weight one less than the weight of $f(x)$.
By induction, $f'(x)$ has weight at least $k=k_0$, and hence
$f(x)$ has weight at least $k_0+1$.

Now suppose that $k\ge p$.
Write $f(x)=(x-1)^{k-k_0}\cdot g(x)$ and
$g(x)=\sum_ig_ix^i$.
Since
$(x-1)^{k-k_0}=(x^p-1)^{(k-k_0)/p}$,
for each $0\le j<p$ we have
\begin{equation*}
\sum_i f_{pi+j}\,x^{pi+j}=
(x-1)^{k-k_0}
\sum_i g_{pi+j}\,x^{pi+j}
\end{equation*}
and, consequently,
\begin{equation}\label{eq:section}
\sum_i f_{pi+j}\,y^{i}=
(y-1)^{(k-k_0)/p}
\sum_i g_{pi+j}\,y^{i},
\end{equation}
having set $y=x^p$.
The polynomial
$\sum_{0\le j<p}\big(\sum_i g_{pi+j}\big)x^{j}$
equals the reduction of $g(x)$ modulo $x^p-1$,
and hence has $1$ as a root with multiplicity exactly $k_0$.
The first part of the proof shows that its weight is at least $k_0+1$.
Hence $\sum_i g_{pi+j}\neq 0$ for at least $k_0+1$ distinct values of $j$.
Consequently, for those values of $j$ the polynomial at the left-hand side of Equation~\eqref{eq:section}
is nonzero, and has $1$ as a root with multiplicity exactly $(k-k_0)/p=\sum_{t>0} k_tp^{t-1}$.
Since this multiplicity is less than $k$, induction implies that
$\sum_i f_{pi+j}\,y^{i}$ has weight at least
$\prod_{t>0}(k_t+1)$.
We conclude that $f(x)$ has weight at least
$(k_0+1)\cdot\prod_{t>0}(k_t+1)=\prod_{t\ge 0}(k_t+1)$.
\end{proof}

Now we return to characteristic zero and give our alternative proof of Lemma~\ref{lemma:weight}.
In order to avoid certain alternating signs, from now on it will be convenient to choose $-1$ instead of $1$
as our distinguished nonzero root of $f(x)$ with multiplicity $k$.

\begin{proof}[Proof of Lemma~\ref{lemma:weight}]
Let $f(x)=\sum_{i=0}^{n}f_ix^i$ be monic of degree $n$, and suppose that
$f(x)=(x+1)^kg(x)$, with $g(x)=\sum_{i=0}^{n-k}g_ix^i$.
Hence we have
\begin{equation}\label{eq:f_i}
f_i=\sum_{j=0}^{n-k}g_j\binom{k}{i-j}
\end{equation}
for all $i=0,\ldots,n$.
Now note that
\begin{align*}
\binom{k}{i-j}&=\frac{k!}{(i-j)!(k+j-i)!}
\\
&=\frac{k!}{i!(n-i)!}\cdot\frac{i!}{(i-j)!}\cdot\frac{(n-i)!}{(k+j-i)!}
\\
&=\frac{k!}{i!(n-i)!}\cdot(i)_j\cdot(n-i)_{n-k-j},
\end{align*}
where $(a)_b$ denotes the integer $a!/(a-b)!=a(a-1)\cdots(a-b+1)$, for $0\le b\le a$.
Therefore, we can rewrite Equation~\eqref{eq:f_i} in the form
\[
i!(n-i)!\cdot f_i=
k!\sum_{j=0}^{n-k}g_j\cdot(i)_j\cdot(n-i)_{n-k-j}.
\]
In particular, $f_i$ vanishes, for some $i$ with $0\le i<n$, if and only if
\[
\sum_{j=0}^{n-k}g_j\cdot(i)_j\cdot(n-i)_{n-k-j}
\]
vanishes.
This expression can be viewed as a polynomial in the indeterminate $i$, of degree at most $n-k$,
and it is not the zero polynomial, otherwise $f_i$ would vanish for all $i<n$.
Hence it has at most $n-k$ roots, and so at most $n-k$ of the coefficients $f_i$ vanish.
Consequently, $f(x)$ has weight at least $k+1$.
\end{proof}

\begin{proof}[Proof of Theorem~\ref{thm:count_extremal}]
Consider an arbitrary polynomial of the form $f(x)=(x+1)^kg(x)$ and of degree at most $n$.
Let $I$ be a subset of $\{0,1,\ldots,n-1\}$.
Continuing the notation of our Proof of Lemma~\ref{lemma:weight},
and according to~\eqref{eq:f_i}, the condition that $f_i=0$ for all $i\in I$
is equivalent to
\begin{equation}\label{eq:system}
\sum_{j=0}^{n-k-1}g_j\binom{k}{i-j}=-g_{n-k}\binom{k}{i-n+k}
\qquad\text{for all $i\in I$}.
\end{equation}
Now suppose that $|I|=n-k$.
Viewing $g_{n-k}$ as assigned,~\eqref{eq:system} is a system of $n-k$ linear equations
in the $n-k$ indeterminates $g_0,\ldots,g_{n-k-1}$.
We claim that the system has a unique solution.

In fact, if $g(x)$ is the polynomial corresponding to any solution of the homogeneous system
obtained by setting $g_{n-k}=0$ in~\eqref{eq:system},
then $f(x)=(x+1)^kg(x)$ has (degree less than $n$ and)
weight at most $k$.
Because of Lemma~\ref{lemma:weight} then $f(x)$ can only be the zero polynomial,
and hence the system has maximal rank $n-k$.
If we now set $g_{n-k}=1$ in~\eqref{eq:system}, then the system has a unique solution, as claimed.
This solution yields a unique monic polynomial $f(x)=(x+1)^kg(x)$ of weight at most $k+1$, and hence
exactly $k+1$ according to Lemma~\ref{lemma:weight}.

According to Lemma~\ref{lemma:weight}, the system of equations~\eqref{eq:system}
can be satisfied for a nonzero polynomial only if $|I|\le n-k$.
Consequently, distinct subsets $I$ of $\{0,1,\ldots,n-1\}$ with $|I|=n-k$ yield distinct
monic polynomials $(x+1)^kg(x)$ of degree $n$ and weight $k+1$.
We conclude that there are exactly $\binom{n}{k}$ such polynomials.
\end{proof}

The Proof of Theorem~\ref{thm:count_extremal},
like the preceding Proof of Lemma~\ref{lemma:weight}, remains valid in prime characteristic $p$
provided $n<p$.
We illustrate the obstacles arising when $n\ge p$ by considering the simplest case $n=p$,
which does not fit the conclusion of Theorem~\ref{thm:count_extremal}.

\begin{example}\label{ex:n=p}
Work over a field of characteristic $p$ and suppose that $0<k<n=p$.
The first two paragraphs of the Proof of Theorem~\ref{thm:count_extremal} apply unchanged, but the last paragraph fails.
In particular, the system of equations~\eqref{eq:system} with $|I|=p-k$ has only the null solution when $g_{p-k}=0$,
and hence has a unique solution when $g_{p-k}=1$.
However, the system may have nonzero solutions even when $|I|>p-k$.
Consequently, to each subset $I$ of $\{0,1,\ldots,p-1\}$ with $|I|=p-k$
there does correspond a unique monic polynomial $f(x)=(x+1)^kg(x)$
of degree $n$ and weight {\em at most} $k+1$, but polynomials $f(x)$ associated with distinct subsets $I$ need not be distinct.
In fact, the polynomial $f(x)=x^p+1$ is associated with all those subsets $I$ not containing $0$.
On the other hand, for each polynomial $f(x)$ associated with a subset $I$ containing $0$,
the polynomial $f(x)/x$ has degree $p-1$ and weight at most $k+1$.
Because of Theorem~\ref{thm:count_extremal}, those polynomials $f(x)/x$ are distinct and in number of $\binom{p-1}{k}$.
We conclude that, in characteristic $p$, there are exactly $\binom{p-1}{k}+1$ monic polynomials
which are multiples of $(x+1)^k$, have degree $p$ and weight at most $k+1$.
More precisely, all of them except $x^p+1$ have weight $k+1$ and are multiples of $x$.
\end{example}

\begin{proof}[Proof of Theorem~\ref{thm:count}]
For each subset $J$ of $\{0,1,\ldots,n-1\}$ let $M_J$ be the number of polynomials in $\F_q[x]$
of the form
\[
x^n+\sum_{j\in J}f_jx^j=f(x)=(x+1)^kg(x)
\]
for some $g(x)\in\F_q[x]$, with $f_j\neq 0$ for all $j\in J$.
We will prove that $M_J=M_{|J|}$ for all $J$.

The number of polynomials of the form described above, but without the requirement that
$f_j\neq 0$ for all $j\in J$, equals
$\sum_{J'\subseteq J}M_{J'}$.
However, those polynomials are in bijective correspondence with the solutions
of the linear system~\eqref{eq:system}, with
$I=\{0,1,\ldots,n-1\}\setminus J$ and $g_{n-k}=1$.
We have seen in the Proof of Theorem~\ref{thm:count_extremal}
that any subsystem of~\eqref{eq:system} consisting of $n-k$ equations
has exactly one solution.
Since the system under consideration here comprises $|I|=n-|J|$ equations, it has no solutions if $|J|< k$,
and has $q^{|J|-k}$ solutions otherwise.
Therefore, we have
\[
\sum_{J'\subseteq J}M_{J'}=
\begin{cases}
q^{|J|-k} &\text{if $|J|\ge k$,}
\\
0 &\text{otherwise.}
\end{cases}
\]
An application of M\"obius inversion
(see~\cite[page 202]{Com}, for example)
yields that
\begin{align*}
M_J&=
\sum_{J'\subseteq J,\ |J'|\ge k}
(-1)^{|J\setminus J'|}\ q^{|J'|-k}
\\&=
\sum_{v\ge k}
(-1)^{w-v-1}\binom{w-1}{v}q^{v-k},
\end{align*}
where in the last passage we have set $|J|=w-1$ and $|J'|=v$.
\end{proof}

Corollary~\ref{cor:count} follows at once from Theorem~\ref{thm:count}.
Note that, of course,
\begin{align*}
\sum_{w}\binom{n}{w-1}M_w
&=
\sum_{v\ge k}q^{v-k}\sum_{w}(-1)^{w-v-1}\binom{n}{w-1}\binom{w-1}{v}
\\&=
\sum_{v\ge k}q^{v-k}\sum_{w}(-1)^{n-v}\delta_{n,v}
=
q^{n-k}
\end{align*}
is the total number of monic polynomials of degree $n$
which are multiples of $(x+1)^k$.

\bibliography{References}

\end{document}